\documentclass[12pt]{article}
\usepackage{amsfonts}
\usepackage{polski}
\usepackage{hyphenat}

\begin{document}
\def\R{\mathbb{R}}
\def\C{\mathbb{C}}
\def\Z{\mathbb{Z}}
\def\N{\mathbb{N}}
\def\Q{\mathbb{Q}}
\def\D{\mathbb{D}}
\def\Sp{{\mathbb{S}}}
\def\T{\mathbb{T}}
\def\hb{\hfil \break}
\def\ni{\noindent}
\def\i{\indent}
\def\a{\alpha}
\def\b{\beta}
\def\e{\epsilon}
\def\d{\delta}
\def\De{\Delta}
\def\g{\gamma}
\def\qq{\qquad}
\def\L{\Lambda}
\def\E{\cal E}
\def\G{\Gamma}
\def\F{\cal F}
\def\K{\cal K}
\def\A{\cal A}
\def\B{\cal B}
\def\M{\cal M}
\def\P{\cal P}
\def\Om{\Omega}
\def\om{\omega}
\def\s{\sigma}
\def\th{\theta}
\def\Th{\Theta}
\def\z{\zeta}
\def\p{\phi}
\def\m{\mu}
\def\n{\nu}
\def\Si{\Sigma}
\def\q{\quad}
\def\qq{\qquad}
\def\half{\frac{1}{2}}
\def\hb{\hfil \break}
\def\half{\frac{1}{2}}
\def\pa{\partial}
\def\r{\rho}
\def\Spd{{{\Sp}}^d}
\begin{center}
{\bf EXTREMES AND REGULAR VARIATION}
\end{center}

\begin{center}
{\bf N. H. BINGHAM and A. J. OSTASZEWSKI}
\end{center}

\begin{center}
{\bf Abstract}
\end{center}
\i We survey the connections between extreme-value theory and regular variation, in one and higher dimensions, from the algebraic point of view of our recent work on Popa groups. \\
\ni {\it Key words} Extreme value, Popa group, generalised Pareto distribution, peaks over thresholds, copula, spectral measure, D-norm, dependence structure, max-stable process, spatio-temporal process \\
\ni {\it Mathematics Subject Classification} 60G70, 62G32  \\

\ni {\bf 1.  One dimension} \\
\i The simplest case is that of an independent and identically distributed (iid) sequence $(X_n)$ with law $F$; write
$$
M_n := \max \{X_1, \cdots, X_n \} \ \hbox{or} \ \vee_1^n X_i.
$$
If there are centring constants $b_n$ and norming constants $a_n$ such that
$$
a_n M_n + b_n \to G \quad \hbox{in law} \quad (n \to \infty)
$$
for some non-degenerate probability distribution $G$, then $G$ is called an {\it extreme-value distribution (EVD)} (or {\it extremal law}), and $F$ belongs to the {\it domain of attraction} of $G$, $F \in D(G)$.  The EVD are also the {\it max-infinitely divisible (max-id)} laws [BalR].  \\
\i We are working here to within an affine transformation (this would change the centring and scaling but preserve the limit), that is, to within {\it type} (location and scale) [Loe I, IV.14].  Modulo type, the limits $G$ (in {\it one} dimension) have a simple parametric description (see e.g. [BinGT, Th. 8.13.1]): \\

\ni {\bf Theorem (Fisher-Tippett theorem)}, [FisT], 1928.  To within type, the extremal laws are exactly the following:
$$
{\Phi}_{\a}, \q (\a > 0); \qq {\Psi}_{\a}, \q (\a > 0); \qq \L,
$$
where the Fr\'echet (${\Phi}_{\a}$), Weibull (${\Psi}_{\a}$) and Gumbel ($\Lambda$) laws are given by
$$
{\Phi}_{\a} := 0 \q (x \leq 0), \q \exp \{ - x^{-{\a}} \} \q (x \geq 0);
$$
$$
{\Psi}_{\a} := \exp \{ - (-x)^{\s} \} \q (x \leq 0), \qq 1 \q (x \geq 0);
$$
$$
\L(x) := \exp \{- e^{-x} \} \q (x \in \R).
$$

\i Particularly for statistical purposes, it is often better to combine these three into one parametric family, the {\it generalized extreme value (GEV)} laws (see e.g. [Col, 3.1.3]).  These have one extremal parameter $\a \in \R$ and two type parameters $\m \in \R$ (location) and $\s > 0$ (scale):
$$
G(x)
:= \exp \Bigl( - \Bigl[ 1 + \a \Bigl( \frac{x - \m}{\s} \Bigr) \Bigr]^{-1/{\a}} \Bigr) \qquad \hbox{where} \ \ \Bigl[ \cdots \Bigr] > 0.           \eqno(GEV)
$$
Here $\a > 0$ corresponds to the Fr\'echet ${\Phi}_{\a}$, $\a = 0$ to the Gumbel $\L$ (using $(1 + x/n)^n \to e^x$ as $n \to \infty$; we interpret $e^x$ as the `$n = \infty$', or `$\a = 0$', case, by the `L' Hospital convention') and $\a < 0$ to the Weibull ${\Psi}_{\a}$.  Taking $\mu = 0$ and $\s = 1$ for simplicity (as we may), this gives the {\it extreme-value distributions}
$$
G_{\a}(x) := \exp (- g_{\a}(x)), \qquad g_{\a}(x) := [ 1 + \a x]_+^{-1/{\a}}. \eqno(EVD)
$$
Here the parameter $\a \in \R$ is called the {\it extreme-value index (EVI)} or {\it extremal index}.  The upper end-point $x_+$ of $F$ is $\infty$ for $\a \geq 0$ (with a power tail for $\a > 0$ and an exponential tail for $\a = 0$); for $\a < 0$ $x_+ = -1/\a$, with a power tail to the left of $x_+$. \\
\i The domains of attraction in the Fr\'echet and Weibull cases, due to Gnedenko [Gne] in 1943 (see e.g. [BinGT, Th. 8.13.2,3]) are simple: writing $\overline{F} := 1 - F$ for the tail of $F$ and $R_{\rho}$ for the class of (positive measurable) functions varying regularly at infinity with index $\rho$, \\
(i) $F \in D({\Phi}_{\a})$ iff $\overline{F} \in R_{- \a}$; \\
(ii) $F \in D({\Psi}_{\a})$ iff $F$ has finite upper end-point $x_+$ and
$\overline{F}(x_+ - 1/.) \in R_{-\a}$. \\
The Gumbel case is more complicated (de Haan [dHaa1,2] in 1970-71, [BinGT, Th. 8.13.4]; cf. [BinGT, Ch. 3, De Haan theory]): \\
(iii) $F \in D(\Lambda)$ iff
$$
\overline{F}(t + x a(t))/\overline{F}(t)
\to g_0(x) := e^{-x} \qquad (x \to \infty),                      \eqno(\ast)
$$
for some {\it auxiliary function} $a > 0$, which may be taken [EmbKM, (3.34)] as
$$
a(t) := \int_t^{x_+} \overline{F}(u) du/\overline{F}(t)
\quad (t < x_+),                                                 \eqno(aux)
$$
and satisfies
$$
a(t + x a(t))/a(t) \to 1 \quad (t \to \infty).                   \eqno(Beu)
$$
Such functions are called {\it Beurling slowly varying} (see e.g. [BinGT, \S 2.11], [BinO1] and the references cited there).  If also $(Beu)$ holds uniformly on compact $x$-sets, $a$ is called {\it self-neglecting}, $a \in SN$  (cf. [BieGST, \S 2.5.2]):
$$
a(t + x a(t))/a(t) \to 1 \quad (x \to \infty) \quad \hbox{(uniformly on compact $x$-sets}).                                                       \eqno(SN)
$$
\i An alternative criterion for $D(\Lambda)$ had been given in 1968 by Marcus and Pinsky [MarP]. \\
\i The three domain-of-attraction conditions may be unified (using the L'Hos-\allowbreak pital convention as above) as follows: $F \in D(G_{\a})$ iff
$$
\overline{F}(t + x a(t))/\overline{F}(t)
\to g_{\a}(x) := (1 + \a x)_+^{-1/\a} \quad (t \to \infty)                   \eqno(\ast \ast)
$$
for some auxiliary function $a$, and then
$$
a(t + x a(t))/a(t) \to 1 + \a x \quad (t \to \infty),           \eqno(\a Beu)
$$
extending the $\a = 0$ case $(Beu)$ above  (see e.g. [BieGST, \S 2.6]). \\
\i For a {\it continuous} Beurling slowly varying $a$, $a \in SN$ and $a(x)=o(x)$ (Bloom's theorem: [BinGT, \S 2.11], [BinO1]). The relation $(\a Beu)$ defines the {\it self-equivarying} functions $a \in SE$ [Ost]; here $a(x)/(1 + \a x) \in SN$ and so allows for $a(x)=O(x)$ (cf. the case $a(x):=1+ \a x$ with $\a > 0$). \\

\ni {\it Von Mises conditions} \\
\i In 1936, von Mises [vMis] gave {\it sufficient} conditions for membership of these domains of attraction, assuming that $F$ has a density $f$ (there is no essential loss of generality here; see below).  We formulate these in terms of the {\it hazard rate} $h$ of survival analysis (see e.g. [CoxO, \S 2.2]):
$$
h(x) := f(x)/\overline{F}(x) = f(x)/\int_x^{\infty} f(u) du.
$$
Below, we shall also need the {\it inverse hazard function}
$$
i(x) := 1/h(x) = \int_x^{\infty} f(u) du/f(x).
$$
Observe that (when the density $f$ exists, as here) the numerator and denominator in $(aux)$ are the integrals of those here.  As one may integrate (though not necessarily differentiate) asymptotic relations, we infer that when $i$ exists it may be used as an auxiliary function $a$ as in $(aux)$. \\
\i Recall the Smooth Variation Theorem [BinGT, \S 1.8]: in any situation in regular variation, (one is working to within asymptotic equivalence
$\sim$, and so) there is no essential loss in assuming that $F$ has a density (even a $C^{\infty}$ density) $f$.  Indeed, Balkema and de Haan [BaldH2] show that in all three cases, if $F \in D(G)$ for $G$ an extremal law, then $\overline{F} \sim \overline{F_*}$, where $F_*$ satisfies a von Mises condition. \\
\i The von Mises conditions in the three cases are (a)-(c) below. \\
(a) For $\Phi_{\a}$: if $x_+ = \infty$ and
$$
x h(x) = x/i(x) \to \a > 0 \quad (x \to \infty),                  \eqno(vM\Phi)
$$
then $F \in D({\Phi}_{\a})$. \\
\i That (a) is equivalent to (i) in the density case follows by Karamata's Theorem [BinGT, \S 1.6]; [BinGT, Th. 8.13.5]. \\
(b) For ${\Psi}_{\a}$: if $x_+ < \infty$ and
$$
(x_+ - x) h(x)
= (x_+ - x)/i(x) \to \a > 0 \quad (x \to \infty),               \eqno(vM\Psi)
$$
then $F \in \Psi(\a)$. \\
\i The proof uses (ii) as above [BinGT, Th. 8.13.6]. \\
(c) Taking $x_+ = \infty$ for simplicity: if
$$
i'(x) \to 0 \quad (x \to \infty),                              \eqno(vM\Lambda)
$$
then $F \in D(\Lambda)$ [BaldH1].  The proof [BinGT, Th. 8.13.7] hinges on [BinGT, Lemma 8.13.8]: if $a(.) > 0$ and $a'(t) \to 0$ as $t \to \infty$, then $a \in SN$.  This actually characterises $SN$: the representation theorem for $a \in SN$ is [BinO1, Th. 9]
$$
a(x) = c (1 + o(1)) \int_0^x e(u) dy, \quad e \in C^1,
\ e(x) \to 0 \ \ (x \to \infty).
$$
\i Rates of convergence in the above were studied by Falk and Marohn [FalM]. \\

\ni {\it Von Mises functions} \\
\i Call a distribution function $F$ a {\it von Mises function} with {\it auxiliary function} $a$ if [EmbKM, \S 3.3.3] for some $c, d \in (0,\infty)$,
$$
\overline{F}(x) = c \ \exp \{ -\int_d^x dt/a(t) \},
\quad a'(t) \to 0 \ \ (t \to \infty).
$$
Then (as above) one can take the auxiliary function $a$ as the inverse hazard function $i$ above, or (see below) the mean excess function $e$ (when it exists).  One can pass to full generality by replacing the constant $c$ above by a function $c(x) \to c$: $D(\Lambda)$ consists of von Mises functions and their tail-equivalent distributions [EmbKM, p.144].  And (from $a' \to 0$): when $x_+ = \infty$, tails in $D(\Lambda)$ decrease faster than any power [EmbKM, p.139].  Example: the standard normal law (take $a = i$ and use Mill's ratio). \\

\ni {\it Peaks over thresholds (POT)} \\
\i As always in extreme-value theory, one has two conflicting dangers.  The maxima --- the very high values --- are rare, and focussing on them discards information and may leave too little data.  But if one over-compensates for this by including too much data, one risks distorting things as the extra data is also informative about the distribution away from the tails.  One approach is to choose a large threshold (which the statistician may choose), $u > 0$ say, and look only at the data exceeding $u$.  These are the {\it peaks over thresholds (POT)}.  Here one focusses on the {\it exceedances} $Y = X - u$ when positive, and their conditional law $F_u$ given $X > u$.  This leads to
\begin{eqnarray*}
F_u(x) := P(Y > u + x a(u)| Y > 0)
&=& P \Bigl( \frac{X - u}{a(u)} > x|X > u \Bigr) \\
&=& \overline{F}(u + x a(u))/\overline{F}(u) \\
&\to & g_{\a}(x) := (1 + \a x)_+^{-1/\a} \quad (u \to \infty),
\end{eqnarray*}
as in $(\ast \ast)$ above.  Thus the conditional distribution of $(X - u)/a(u)|X > u$ has limit
$$
H_{\a}(x) := 1 + \log G_{\a}(x)
= 1 - g_{\a}(x) = 1 - (1 + \a x)_+^{-1/\a},          \eqno(GPD)
$$
the {\it generalised Pareto distribution (GPD)} (`EVD for max, GPD for POT'). \\
\i There are several ways of motivating the use of GPD: \\
(i) Pickands [Pic1] showed in 1975 that $F_u$ has GPD $H_{\a}$ as limit law iff it has the corresponding EVD as limit of its maxima, i.e. $F \in D(G_{\a})$.  This (in view of [BaldH2]) is the {\it Pickands-Balkema-de Haan theorem} [McNFE, Th. 7.20].\\
(ii) There is {\it threshold stability}: if $Y$ is GP and $u > 0$, then the conditional law of $Y - u|Y > u$ is also GP, and this characterises the GPD.  This property is useful in applications; see e.g. [McNFE, \S 7.2.2]. \\
(iii) If $N$ is Poisson and $(Y_1, \cdots, Y_N)|N$ are iid GP, then $\max (Y_1, \cdots, Y_N)$ has the corresponding EVD; again, this characterises GPD. \\
For details, see e.g. Davison and Smith [DavS \S2] and the references there. \\
\i Statistical work in the one-dimensional setting here centres on the estimation of the extreme-value index $\a$.  One of the commonest estimators here is {\it Hill's estimator} [Hil]; see e.g. [McNFE \S 7.2.4], [BieGST, \S 9.5.2]. \\

\ni {\it Mean excess function} \\
\i When the mean of $X$ exists, the {\it mean excess} (or {\it mean exceedance}) {\it function} of $X$ over the threshold $u$ exists and is
$$
e(u) := E[X - u | X > u].
$$
Integrating by parts,
$$
e(u) = \int_u^{\infty} (x-u)dF(x)/\overline{F}(u)
= - \int_u^{\infty} (x-u)d \overline{F}(x)/\overline{F}(u)
= \int_u^{\infty} \overline{F}(x) dx/\overline{F}(u),
$$
which by $(aux)$ is the general form of the auxiliary function $a$.  Thus, {\it when} $e$ {\it exists, one may take it as the auxiliary function} $a$ (in preference to the inverse hazard function $i$, if preferred). \\

\ni {\it Self-exciting processes} \\
\i One way to relax the independence assumption is to allow {\it self-exciting processes}, where an occurrence makes other occurrences more likely.  This is motivated by aftershocks of earthquakes, but also relevant to financial crises.  This uses {\it Hawkes processes} [Haw]; see [McNFE \S 7.4.3], and for point-process background, [Res1]. \\

\ni {\it Popa groups} \\
\i Referring to $(\ast)$ and $(\a Beu)$, these can now be recognised as the relevant instance of a {\it Popa group} (for Popa groups and general regular variation, see [BinO2]).  The signature is the argument $t + xa(t)$, where the auxiliary function $a$ is self-neglecting.  See the $3 \times 3$ table in [BinO3, Th. BO] (relevant here is the top right-hand corner with $\kappa = -1$).  Likewise, the limit in $(\ast \ast)$ gives the (2,3) (or middle right) entry in the table, with $\kappa = -1$, and after taking logs, the (2,1) (or middle left) entry:
$$
\log \overline{F}(t + x a(t)) - \log \overline{F}(t)
\to {-1/\a} \log (1 + \a x)_+ \quad (t \to \infty),
                                                         \eqno(\ast \ast \ast)
$$
exactly of the form studied in [BinO2] (there the RHS is called the {\it kernel}, $K(x)$).  This leads ([BinO2, \S 2], with $\alpha$ for the $\rho$ there and ${\circ}_{\alpha}$ for the Popa operation) to a {\it Goldie equation} [BinO2, \S 5]
$$
K(x {\circ}_{\alpha} y) = K(x) +g_{\alpha}
(x)K(y);
$$
here
$$\qquad g_{\alpha}(x {\circ}_{\alpha} y) =g_{\alpha}(x)g_{\alpha}(y).$$
\i The authors in [EmbKM] remark (e.g. their p.140) that regular variation `does not seem to be the right tool' for describing von Mises functions.  The general regular variation of [BinO1,2] does seem to be the right tool here, including as it does the Karamata, Bojanic-Karamata/de Haan and Beurling theories of regular variation. \\

\ni {\bf 2. Higher dimensions} \\
\i For general references for multidimensional EVT, see e.g. [EmbKM, Ch. 3,5,6], [BieGST, Ch. 8], [dHaaF, Part II]. [dHaaRe], [FalHR], [Fal1,2].  For multidimensional regular variation, see e.g. Basrak et al. [BasDM]. \\
\i The situation in dimensions $d > 1$ is different from and more complicated than that for $d = 1$ above, regarding both EVT and Popa theory.  We hope to return to such matters elsewhere. \\

\ni {\it Copulas} \\
\i  The theory above extends directly from dimensions 1 to $d$, if each $X$ or $x$ above is now interpreted as a $d$-vector.   Then, as usual in dimension $d > 1$, we split the $d$-dimensional joint distribution function $F$ into the marginals $F_1, \cdots, F_d$, and the {\it copula} $C$ (a probability law on the $d$-cube $[0,1]^d$ with uniform marginals on $[0,1]$), which encodes the dependence structure via {\sl Sklar's theorem} ([Skl]; see e.g. [McNFE, Ch. 5]):
$$
F(x_1, \cdots, x_d) = C(F_1(x_1), \cdots, F_d(x_d)).              \eqno(Skla)
$$
In particular, this shows that one may standardise the marginals $F_i$ in any convenient way, changing only the joint law $F$ but not the copula (dependence structure).  One choice often made in extreme-value theory is to transform to {\it standard Fr\'echet} marginals,
$$
F_i(x) = \exp \{ -1/x \} \quad (x > 0).
$$
When this is done, the EV law is called {\it simple} [BieGST, \S 8.2.2]. \\
\i The limit distributions that can arise are now the {\it multivariate extreme-value (MEV)} laws.  Their copulas link MEV laws with their GEV margins.  With ${\bf u}^t := (u_1^t, \cdots, u_d^t)$ for $t > 0$, these, the {\it EV copulas}, denoted by $C_0$, are characterised by their {\it scaling relation} [McNFE, Th. 7.44]
$$
C_0({\bf u}^t) = C_0^t({\bf u}) \quad (t > 0).                        \eqno(Sca)
$$
\i By analogy with stable laws for sums, a law $G$ is {\it max-stable} if
$$
G^n(a_n x + b_n) \equiv G(x) \quad (n \in \N)
$$
for suitable centring and scaling sequences $(b_n), (a_n)$; these are the GEV laws. \\

\ni {\it Survival copulas} \\
\i In extreme-value theory, it is the upper tails that count.  Taking operations on vectors componentwise and writing
$$
F(x) := P(X \leq x), \qquad \overline{F}(x) := P(X < x),
$$
one can rewrite Sklar's theorem in terms of {\it survival functions}
$\overline{F}$, $\overline{F_i}$: a $d$-dimensional survival function $\overline{F}$ has a {\it survival copula} $C$ with [McNN, Th. 2.1]
$$
\overline{F}(x) = C(\overline{F_1}(x_1), \cdots, \overline{F_d}(x_d)), \qquad
C(u) = \overline{F}(\overline{F_1}^{-1}(u_1), \cdots \overline{F_d}^{-1}(u_d))
$$
($F, F_i$ are 1 at $+ \infty$; $\overline{F}, \overline{F_i}$ are 0 at $+ \infty$; they are accordingly often studied for $x \leq 0$, $x \to - \infty$ rather than $x \geq 0$, $x \to \infty$). \\

\ni {\it Copula convergence} \\
\i The question of multivariate domains of attraction (MD, or MDA) decomposes into those for the marginals and for the copula by the {\it Deheuvels-Galambos theorem} ([Deh], [Gal]; [McNFE, Th. 7.48]): with $F$ as above, $F \in MD(H)$ with
$$
H(x_1, \cdots, x_d) := C_0(H_1(x_1), \cdots, H_d(x_d)),
$$
an MEV law with GEV marginals $H_i$ and EV copula $C_0$ iff \\
(i) $F_i \in MD(H_i)$, $i = 1, \cdots, d$; \\
(ii) $C \in CD(C_0)$, i.e.
$$
C^t(u_1^{1/t}, \cdots, u_d^{1/t})
\to C_0(u_1, \cdots, u_d) = C_0({\bf u}) \quad (t \to \infty)
\quad ({\bf u} \in [0,1]^d).
$$

\ni {\it Peaks over thresholds (POT)} \\
\i The first (and most important) two of the three properties above of POT in one dimension extend to $d$ dimensions.  The first is the $d$-dimensional version of the Pickands-Balkema-de Haan theorem, linking EVD and GPD; the second is threshold stability.  For details, see Rootz\'en and Tajvidi [RooT], Rootz\'en, Segers and Wadsworth [RooSW1,2], Kiriliouk et al. [KirRSW]. \\

\ni {\it Spectral representation} \\
\i The scaling property $(Sca)$ suggests using spherical polar coordinates, $x = (r,\th)$ say ($x \in {\R}_+^d$, $r > 0$, $\th \in {\Sp}_+^{d-1}$).  Then the EV law $G$ has (with $\wedge$ for $\min$) a {\it spectral representation}
$$
\log G(x)
= \int_{{\Sp}_+} {\wedge}_{i=1}^d \Bigl(\frac{{\th}_i}{\Vert \th \Vert} \log G_i(x_i) \Bigr) \ dS(\th)
\quad (x = (r, \th) \in {\R}^d),
$$
where the {\it spectral measure} $S$ satisfies
$$
\int_{{\Sp}_+} \frac{{\th}_i}{\Vert \th \Vert} dS(\th) = 1 \quad (i = 1, \cdots, d)
$$
(see e.g. [dHaa3], [BieGST, \S 8.2.4], [MaoH]).  The regular-variation (or other limiting) properties are handled by the radial component, the dependence structure by the spectral measure. \\

\ni $D$-{\it norms} \\
\i The standard (i.e. with unit Fr\'echet marginals) max-stable (SMS) laws are those with survival functions of the form
$$
\exp \{ - \Vert x \Vert \} \quad (x \leq 0 \in {\R}^d),
$$
for some norm, called a $D$-{\it norm} (`D for dependence', as this norm encodes the dependence structure).  For a textbook treatment, see Falk [Fal2]. \\

\ni {\it Pickands dependence function} \\
\i An EV copula may be specified by using the {\it Pickands dependence function}, $B$ ([Pic2]; [McNFE, Th. 7.45]): $C$ is a $d$-dimensional EV copula iff it has the representation
$$
C({\bf u}) = \exp \Bigl( B \Bigl(\frac{\log u_1}{\sum_1^d u_i}, \cdots,
\frac{\log u_d}{\sum_1^d u_i} \Bigr) \sum_1^d \log u_i \Bigr),
$$
where with $S_d$ the $d$-simplex $\{ x: x_i \geq 0, \sum_1^d x_i = 1 \}$,
$$
B(w) = \int_{S_d} \max (x_1 w_1, \cdots, x_d w_d) dH(x)
$$
with $H$ a finite measure on $S_d$.  Of course, the $d$-simplex needs only $d-1$ coordinates to specify it; this simplification is most worthwhile when $d = 2$ (below). \\

\ni {\it  Two dimensions} \\
\i Things can be made more explicit in two dimensions.  For the theory above, one obtains the representation
$$
C(u_1,u_2) = \exp \{ (\log u_1 + \log u_2) A \bigl( \frac{\log u_1}{\log u_1 + \log u_2} \bigr) \},
$$
where
$$
A(w) = \int_0^1 max ((1-x)w, x(1-w)) dH(x),
$$
for $H$ a measure on $[0,1]$.  The Pickands dependence function $A$ here is characterised by the bounds
$$
\max (w, 1-w) \leq A(w) \leq 1 \quad (0 \leq w \leq 1)
$$
and being (differentiable and) convex. \\

\ni {\it Archimedean copulas} \\
\i In $d$ dimensions, the {\it Archimedean copula} $C$ with {\it generator} $\psi$ is given by
$$
C(u) = \psi({\psi}^{-1}(u_1) + \cdots + {\psi}^{-1}(u_d)).
$$
Here [McNN] $\psi(0) = 1$, $\psi(x) \to 0$ as $x \to \infty$ and $\psi$ is $d$-monotone (has $d-2$ derivatives alternating in sign with $(-)^{d-2} {\psi}^{(d-2)}$ nonincreasing and convex); this characterises Archimedean copulas.  In particular, for $d = 2$,
$$
C(u,v) = \psi({\psi}^{-1}(u) + {\psi}^{-1}(v))
$$
is a copula iff $\psi$ is convex.  \\
\i An alternative to spherical polars uses the $d$-simplex $S_d$ in place of ${\Sp}_+^d$.  This leads to ${\ell}_1$-{\it norm symmetric distributions}, or {\it simplex distributions}, and the {\it Williamson transform}; see [McNN, \S 3].  Here the Archimedean generator $\psi$ is the Williamson transform of the law of the radial part $R$, and one can read off the domain-of-attraction behaviour of $R$ from regular-variation conditions on $\psi$ [LarN, Th. 1].  The dependence structure is now handled by the {\it simplex measure} [GenNR].\\
\i This feature that one `radial' variable handles the tail behaviour and regular-variation aspects, while the others handle the dependence structure, has led to `one-component regular variation' in this context; see Hitz and Evans [HitE]. \\
\i The marginals may require different normalisations; for background here, see e.g. [Res2, \S 6.5.6] (`standard v. non-standard regular variation'). \\

\ni {\it Gumbel copulas} \\
\i For $\th \in [1,\infty)$, the {\it Gumbel copula} with parameter $\th$ is the Archimedean copula with generator $\psi(x) = \exp \{-x^{1/\th} \}$.  This is the only copula which is both Archimedean and extreme-value ([GinR]; [McNN, Cor. 1]). \\

\ni {\it Archimax copulas} \\
\i Call $\ell$ a ($d$-variate) {\it stable tail-dependence function} if for $x_i \geq 0$,
$$
\ell(x_1, \cdots, x_d) = - \log C_0(e^{-x_1}, \cdots, e^{-x_d})
$$
for some extreme-value copula $C_0$.  The Archimax copulas [CharFGN] are those of the form
$$
C_{\psi, \ell}(u_1, \cdots, u_d)
:= \psi \circ \ell ({\psi}^{-1}(u_1), \cdots, {\psi}^{-1}(u_d)).
$$
This construction does indeed yield a copula [CharFGN], and in the case $d = 2$ gives the Archimedean copulas ($A(.) \equiv 1$ above) and the extreme-value copulas ($\psi(t) = e^{-t}$), whence the name. \\

\ni {\it Dependence structure} \\
\i Particularly when $d$ is large, the spectral measure above may be too general to be useful in practice, and so special types of model are often used, the commonest being those of Archimedean type.  While convenient, Archimedean copulas are {\it exchangeable}, which of course is often not the case in practice (`sea and wind').  The arguments of the copula typically represent covariates, and these are often related by conditional independence relationships; these may be represented graphically (see e.g. the monograph by Lauritzen [Lau], and for applications to extremes, Engelke and Hitz [EngH]; see also [HitE]).  Hierarchical relationships between the covariates (`phylogenetic trees') may be represented by hierarchical Archimedean copulas; see e.g. Cossette et al. [CosGMR].  Special types of graphs (vines) occur in such contexts; see e.g. Chang and Joe [ChaJ], Joe et al. [JoeLN], Lee and Joe [LeeJ]. \\

\ni {\it Max-stable processes} \\
\i The case of infinitely many dimensions -- stochastic processes -- is just as important as the case $d < \infty$ above (the classic setting here is the whole of the Dutch coastline, rather than just coastal monitoring stations).  For theory here, see e.g. [dHaaF, Part III]. \\
\i A process $Y$ is {\it max-stable} if when $Y_i$ are independent copies of $Y$, \\
$\max \{Y_1, \cdots, Y_r \}$ has the same distribution as $rY$ for each $r \in {\N}_+$.  These have a spectral representation, for which see de Haan [dHaa3].  For estimation of max-stable processes, see e.g. Chan and So [ChanS]. \\

\ni {\it Spatio-temporal processes} \\
\i Spectral representations have useful interpretations for modelling spatio-temporal processes, e.g. for the {\it storm-profile process} or {\it Smith process} ([Smi1]; [ChanS]):
$$
Z(x) = {\max}_i \ \phi(x - X_i){\Gamma}_i,
$$
where $Z(x)$ represents the maximum effect at location $x$ over an infinite number of storms centred at random points $X_i$ (forming a homogeneous Poisson process on ${\R}^d$) of strengths ${\Gamma}_i$ (a Poisson point process of rate 1), the effect of each being $\phi(t - X_i){\Gamma}_i$ (here $\phi$ is a Gaussian density function with mean 0 and covariance matrix $\Sigma$, whose contours represent the decreasing effect of a storm away from its centre).  Thus the process measures `the all-time worst (effect), here'.  Perhaps such models could be used to describe e.g. the bush fires currently threatening Australia.\\
\i Spatio-temporal max-stable processes in which the space-time spectral function decouples into ones for time and for space given time are given by Embrechts, Koch and Robert [EmbKR] (see also [KocR]).  This allows for the different roles of time and space given time to be modelled separately.  They also allow space to be a sphere, necessary for realistic modelling on a global scale. \\

\ni {\it Tail dependence} \\
\i Studying asymptotic dependence in multivariate tails is important in, e.g., risk management, where one may look to diversify by introducing negative correlation.  For very thin tails (e.g., Gaussian) this is not possible in view of {\it asymptotic independence} (Sibuya [Sib]).  But with heavier tails, tail dependence coefficients are useful here; see e.g. [SchmS], [LarN, \S 5]. \\

\ni {\it Applications} \\
\i For more on spatial processes (random fields), spatio-temporal processes and applications to such things as weather, see e.g. Smith [Smi1], Schlather [Schl], Cooley et al. [CooNN], Davison et al. [DavPR], Davis et al. [DavKS1,2], Sharkey and Winter [ShaW], Abu-Awwad et al. [AbuMR]. \\
\i An extended study of sea and wind, applied to the North Sea flood defences of the Netherlands, is in de Haan and de Ronde [dHaaRo]. \\
\i Particularly with river networks, the spatial relationships between the points at which the data is sampled is crucial.  For a detailed study here, see Asadi et al. [AsaDE].  \\
\i  For financial applications (comparison of two exchange rates), see [McNFE, Ex. 7.53].     \\

\ni {\it Statistics} \\
\i The great difference between one and higher dimensions in the statistics of extreme-value theory is that in the former, {\it parametric} methods suffice (whether one works with EVD or with GPD).  In the latter, one has $d$ such one-dimensional parametric problems (or one $d$-dimensional one) for the marginals, and a non-parametric one for the copula.  The problem is thus {\it semi-parametric}, and may be treated as such (cf. [KluKP], [BicKRW]).  But our focus here is on the copula, which needs to be estimated nonparametrically; see e.g. [DavS], [GudS1,2], [dFonD] (cf. [EasHT], [PapT]), and in two dimensions, [Seg], [GuiPS].  \\
\i For peaks over thresholds in higher dimensions, see e.g. [KirRSW]; for graphical methods, see [LeeJ]. \\

\ni {\bf 3.  Historical comments} \\
1.  The extremal laws are known as the {\it Fr\'echet} (heavy-tailed, ${\Phi}_{\a}$), {\it Gumbel} (light-tailed, $\L$) and {\it Weibull} (bounded tail, ${\Psi}_{\a}$) distributions, after Maurice Fr\'echet (1878-1973), French mathematician, in 1937, Emil Julius Gumbel (1891-1966), German statistician, in 1935 and 1958, and Waloddi Weibull (1887-1979), Swedish engineer, in 1939 and 1951.\hb
2. The Pareto distributions are named after Vilfredo Pareto (1848-1923), Italian economist, in 1896. \\
3. The remarkable pioneering work of Fisher and Tippett [FisT] in 1928 of course pre-dated regular variation, which stems from Karamata in 1930.  \\
4. The remarkable pioneering work of von Mises [vMis] in 1936 did not use regular variation, perhaps because he was not familiar with the journal Karamata published in, {\sl Mathematica (Cluj)}; perhaps because what Karamata was then famous for was his other 1930 paper, on the (Hardy-Littlewood-)Karamata Tauberian theorem for Laplace transforms --- analysis, while von Mises was an applied mathematician. \\
5. The pioneering work of Gnedenko [Gne] in 1943 on limits of maxima also did not use regular variation; nor did the classic monograph of Gnedenko and Kolmogorov [GneK] of 1949.  As a result, the analytic aspects of both were excessively lengthy, tending to mask the essential probabilistic content. \\
6. The subject of extreme-value theory was made much more important by the tragic events of the night of 31 January - 1 February 1953.  There was great loss of life in the UK, and much greater loss in the low-lying Netherlands (see e.g. [dHaa4]). \\
7. The realisation that regular variation was the natural language for limit theorems in probability is due to Sudakov in 1955 (in Volume 1 of {\sl Theory of Probability and its Applications}).  But this was not picked up at the time, and was rediscovered by Feller in Volume II of his book (1966 and 1971).  See e.g. [Bin2] for details. \\
8. Beurling slow variation appeared in Beurling's unpublished work of 1957 on his Tauberian theorem.  See e.g. [Bin1], [BinO1, \S 10.1] for details and references. \\
9. The first systematic application of regular variation to limit theorems in probability was de Haan's 1970 thesis [dHaa1].  This has been the thread running through his extensive and influential work for the last half-century.\\
10. The Balkema-de Haan paper [BaldH2] of 1974 was explicitly a study of applications of regular variation, and in `great age' set the stage for `high thresholds'. \\
11. Threshold methods were developed by hydrologists in the 1970s.  Their theoretical justification stems from Pickands's result ([Pic1], Pickands-Balkema-de Haan theorem), giving a sense in which $F_u$ is well-approximated by some GPD iff $F$ lies in the domain of attraction of some EVD (cf. [Smi2, \S 3]).  Full references up to 1990 are in [DavS].\\

\ni {\bf Postscript} \\
\i It is a pleasure to contribute to this volume, celebrating Ron Doney's 80th birthday.  His long and productive career in probability theory has mainly focused on random walks and (later) L\'evy processes, essentially the limit theory of sums.  Extreme-value theory is essentially the limit theory of maxima.  Sums and maxima have many points of contact (see e.g. [BinGT, \S 8.15]), recently augmented by the fine paper by Caravenna and Doney [CarD] related to the Garsia-Lamperti problem (see e.g. [BinGT, \S \S 8.6.3, 8.7.1]).  A personal point of contact with extremes came for Ron with flooding and the partial collapse on 6 August 2019 of the dam at Whaley Bridge near his home.  Many residents had to be evacuated, fortunately not including Ron and Margaret.  Any such incident stands as a riposte to climate-change deniers everywhere.  Of course, Australia is much in our minds at the time of writing. \\
\i We thank the editors for their kind invitation to contribute to this Festschrift.  The first author thanks the organisers of Extreme Value Analysis 11 for their invitation to speak at EVA11 in Zagreb in July 2019. \\
\i Both authors send their very best wishes to Ron and Margaret. \\

\begin{center}
{\bf References}
\end{center}
\ni [AbuMR] A.-F. Abu-Awwad, V. Maume-Deschampts and P. Ribereau, Semi-parametric estimation for space-time max-stable processes: E-madogram-based estimation approach.  arXiv:1905:07912. \\
\ni [AsaDE] P. Asadi, A. C. Davison and S. Engelke, Extremes on river networks.  {\sl Ann. Appl. Stat.} {\bf 9} (2015), 2023-2050. \\
\ni [BaldH1] A. A. Balkema and L. de Haan, On R. von Mises' condition for the domain of attraction of $\exp(- e^{-x})$.  {\sl Ann. Math. Statist.} {\bf 43} (1972), 1352-1354. \\
\ni [BaldH2]  A. A. Balkema and L. de Haan, Residual lifetime at great age.  {\sl Ann. Prob.} {\bf 2} (1974), 792-804. \\
\ni [BalR] A. A. Balkema and S. I. Resnick, Max-infinite divisibility.  {\sl J. Appl. Prob.} {\bf 14} (1977), 309-319. \\
\ni [BasDM] B. Basrak, R. A. Davis and T. Mikosch, A characterization of multidimensional regular variation.  {\sl Ann. Appl. Prob.} {\bf 12}(2012), 908-920. \\
\ni [BicKRW] P. J. Bickel, C. A. J. Klaassen, Y. Ritov and J. A. Wellner, {\sl Efficient and adaptive estimation for semiparametric models}.  The Johns Hopkins University Press, 1993 (Springer, 1998). \\
\ni [BieGST] J. Bierlant, Y. Goegebeur, J. Segers and J. Teugels, {\sl Statistics of extremes: Theory and applications}.  Wiley, 2004. \\
\ni [Bin1] N. H. Bingham, Tauberian theorems and the central limit theorem.  {\sl Ann. Prob.} {\bf 9} (1981), 221-231. \\
\ni [Bin2] N. H. Bingham, Regular variation and probability: The early years.  {\sl J. Computational and Applied Mathematics} {\bf 200} (2007), 357-363 (J. L. Teugels Festschrift). \\
\ni [BinGT] N. H. Bingham, C. M. Goldie and J. L. Teugels, {\sl Regular variation}.  Cambridge University Press, 1987/1989. \\
\ni [BinO1] N. H. Bingham and A. J. Ostaszewski, Beurling slow and regular variation.  {\sl Transactions London Math. Soc.} {\bf 1} (2014), 29-56. \\
\ni [BinO2] N. H. Bingham and A. J. Ostaszewski, General regular variation, Popa groups and quantifier weakening. {\sl J. Math. Anal. Appl.} {\bf 483} (2020), 123610; arXiv:1901.05996. \\
\ni [CarD] F. Caravenna and R. A. Doney, Local large deviations and the strong renewal theorem.  {\sl Electronic J. Prob.} {\bf 24} (2019), 1-48. \\
\ni [ChaJ] B. Chang and H. Joe, Predictions based on conditional distributions of vine copulas.  {\sl Comput. Stat. Data Anal.} {\bf 139} (2019), 45-63.  \\
\ni [ChanS] R. K. S. Chan and M. K. P. So, On the performance of the Bayes composite likelihood estimation of max-stable processes.  {\sl J. Stat. Comp. Simul.} {\bf 87} (2017), 2869-2881. \\
\ni [CharFGN] A. Charpentier, A.-L. Foug\`eres, C. Genest and J. G. Ne\u slehov\'a, Multivariate Archimax copulas.  {\sl J. Multivariate Anal.} {\bf 126} (2014), 118-136. \\
\ni [Col] S. Coles, {\sl An introduction to statistical modeling of extreme values}.  Springer, 2001. \\
\ni [CooNN]  D. Cooley, D. Nychka and P. Naveau, Bayesian spatial modelling of extreme precipitation return levels.  {\sl J. Amer. Stat. Assoc.} {\bf 102} (2007), 824-860.  \\
\ni [CooT] D. Cooley and E. Thibaud, Decompositions of dependence for high-dimensional extremes.  {\sl Biometrika} {\bf 106} (2019), 587-604. \\
\ni [CosGMR] H. Cossette, S.-P. Gadoury, E. Marceau and C. Y. Robert, Composite likelihood estimation methods for hierarchical Archimedean copulas defined with multivariate compound distributions.  {\sl J. Multivariate Anal.} {\bf 172} (2019), 59-83.  \\
\ni [CoxO] D. R. Cox and D. Oakes, {\sl Analysis of survival data}.  Chapman \& Hall/CRC, 1984. \\
\ni [DavKS1] R. A. Davis, C. Kl\"uppelberg and C. Steinkohl, Max-stable processes for modelling extremes observed in space and time.  {\sl J. Korean Stat. Soc.} {\bf 42} (2013), 399-413.  \\
\ni [DavKS2] R. A. Davis, C. Kl\"uppelberg and C. Steinkohl, Statistial inference for max-stable processes in space and time.  {\sl J. Roy. Stat. Soc. B} {\bf 75} (2013), 791-819. \\
\ni [DavPR] A. C. Davison, S. A. Padoan and M. Ribatet, Statistical modelling of spatial extremes (with discussion).  {\sl Stat. Science} {\bf 27} (2012), 161-201. \\
\ni [DavS] A. C. Davison and R. L. Smith, Models for exceedances over high thresholds (with discussion).  {\sl J. Roy. Stat. Soc. B} {\bf 52} (1990), 393-442.\\
\ni [Deh] P. Deheuvels, Characterisation compl\`ete des lois extr\^emes multivarie\'es et la convergence aux types extr\^emes.  {\sl Publ. Inst. Sci. U. Paris} {\bf 23} (1978), 1-36. \\
\ni [EasHT] E. F. Easthoe, J. E. Heffernan and J. A. Tawn, Nonparametric estimation of the spectral measure, and associated dependence measures, for multivariate extreme values using a limiting conditional representation.  {\sl Extremes} {\bf 17} (2014), 25-43. \\
\ni [EmbKM] P. Embrechts, C. Kl\"upelberg and T. Mikosch, {\sl Modelling extremal events for insurance and finance}.  Springer, 1997. \\
\ni [EmbKR] P. Embrechts, E. Koch and C. Robert, Space-time max-stable models with spectral separability.  {\sl Adv. Appl. Prob.} {\bf 48A} (2016) (N. H. Bingham Festschrift), 77-97. \\
\ni [EngH] S. Engelke and A. S. Hitz, Graphical models for extremes.  arXiv:1812: 01734. \\
\ni [Fal1] M. Falk, It was 30 years ago today when Laurens de Haan went the multivariate way.  {\sl Extremes} {\bf 11} (2008), 55-80. \\
\ni [Fal2] M. Falk, {\sl Multivariate extreme-value theory and} $D$-{\sl norms}.  Springer, 2019. \\
\ni [FalHR] M. Falk, J. H\"usler and R.-D. Reiss, {\sl Laws of small numbers: Extremes and rare events}, 3rd ed., Birkh\"auser, 2010 (2nd ed., 2004, 1st ed. 1994). \\
\ni [FalM] M. Falk and F. Marohn, Von Mises conditions revisited.  {\sl Ann. Probab.} {\bf 21} (1993), 1310-1328. \\
\ni [FalPW] M. Falk, S. A. Padoan and F. Wisheckel, Generalized Pareto copulas: A key to multivariate extremes.  {\sl J. Multivariate Anal.} {\bf 174} (2019), 104538. \\
\ni [FisT] R. A. Fisher and L. H. C. Tippett, Limiting forms of the frequency of the largest or smallest member of a sample.  {\sl Proc. Cambridge Phil. Soc.} {\bf 24} (1928), 180-190 (reprinted in {\sl Collected papers of R. A. Fisher} {\bf 2} (1972), 208-219, U. Adelaide Press, 1972). \\
\ni [dFonD] R. de Fondeville and A. C. Davison, High-dimensional peaks-over-threshold inference.  {\sl Biometrika} {\bf 105} (2018), 575-592. \\
\ni [Gal] J. Galambos,  {\sl The asymptotic theory of extreme order statistics}.  Wiley, 1978 (2nd ed., Krieger, 1987). \\
\ni [GenNR] C. Genest, J. G. Ne\u slehov\'a and L.-P. Rivest, The class of multivariate max-id copulas with ${\ell}_1$-norm symmetric exponent measure.  {\sl Bernoulli} {\bf 24}(4B) (2018), 3751-3790. \\
\ni [GenR] C. Genest and L. Rivest, On the multivariate probability integral transform.  {\sl Stat. Prob. Letters} {\bf 53} (2001), 391-399. \\
\ni [Gne] B. V. Gnedenko, Sur la distribution limite du terme maximum d'une s\'erie al\'eatoire.  {\sl Ann. Math.} {\bf 44} (1943), 423-453. \\
\ni [GneK] B. V. Gnedenko and A. N. Kolmogorov, {\sl Limit distributions for sums of independent random variables}.  Addison-Wesley, 1954 (Russian edition, Gostelkhizdat, Moscow, 1949). \\
\ni [GudS1] G. Gudendorf and J. Segers, Extreme-value copulas. {\sl Copula theory and its applications} (Warsaw, 2009, ed. P. Jaworski, F. Durante, W. H\"ardle and W. Rychlik), {\sl Lecture Notes in Statistics} (2010), 127–-146. \\
\ni [GudS2] G. Gudendorf and J. Segers, Nonparametric estimation of an extreme-value copula in arbitrary dimensions, {\sl J. Multivariate Analysis} {\bf 102} (2011), 37–-47. \\
\ni [GuiPS] S. Guillotte, F. Perron and J. Segers, Non-parametric Bayesian inference on bivariate extremes, {\sl J. Royal Statistical Soc. B} {\bf 73} (2011), 377–-406. \\
\ni [dHaa1] L. de Haan, {\sl On regular variation and its application to the weak convergence of sample extremes}.  Math. Centre Tract {\bf 32}, Amsterdam, 1970.  \\
\ni [dHaa2] L. de Haan, A form of regular variation and its application to the domain of attraction of the double exponential distribution.  {\sl Z. Wahrschein.} {\bf 17} (1971), 241-258.  \\
\ni [dHaa3] L. de Haan, A spectral representation of max-stable processes.  {\sl Ann. Prob.} {\bf 12} (1984), 1194-1204. \\
\ni [dHaa4] L. de Haan, Fighting the arch-enemy with mathematics.  {\sl Statistica Neerlandica} {\bf 44} (1990), 45-68. \\
\ni [dHaaF] L. de Haan and A. Ferreira, {\sl Extreme value theory: An introduction}.  Springer, 2006. \\
\ni [dHaaRe] L. de Haan and S. I. Resnick, Limit theory for multivariate sample extremes.  {\sl Z. Wahrschein.} {\bf 40} (1977), 317-337. \\
\ni [dHaaRo] L. de Haan and J. de Ronde, Sea and wind: multivariate extremes at work.  {\sl Extremes} {\bf 1} (1998), 7-45. \\
\ni [Haw] A. G. Hawkes, Point spectra of some self-exciting point processes.  {\sl  J. Roy. Statist. Soc. B} {\bf 33} (1971), 438-443. \\
\ni [Hil] B. M. Hill, A simple general approach to inference about the tail of a distribution.  {\sl Ann. Statist.} {\bf 3} (1975), 1163-1174. \\
\ni [HitE] A. Hitz and R. Evans, One-component regular variation and graphical modelling of extremes.  {\sl J. Appl. Prob.} {\bf 53} (2016), 733-746. \\
\ni [JoeLN] H. Joe, H. Li and A. K. Nikoloulopoulos, Tail dependence functions and vine copulas.  {\sl J. Multivariate Anal.} {\bf 101} (2010), 252-270.  \\
\ni [KirRWS] A. Kiriliouk, H. Rootz\'en, J. L. Wadsworth and J. Segers, Peaks-over-thresholds modeling with multivariate generalized Pareto distributions, {\sl Technometrics} {\bf 61} (2019), 123–-135. \\
\ni [KluKP] C. Kl\"uppelberg, G. Kuhn and L. Peng, Semi-parametric models for the multivariate tail-dependence function -- the asymptotically dependent case.  {\sl Scand. J. Stat.} {\bf 35} (2008), 701-718. \\
\ni [KocR] E. Koch and C. Y. Robert, Geometric ergodicity for some space-time max-stable Markov chains.  {\sl Stat. Prob. Letters} {\bf 145} (2019), 43-49.\\
\ni [LarN] M. Larsson and J. N\u eslehov\'a, Extremal behaviour of Archimedean copulas.  {\sl Adv. Appl. Prob.} {\bf 43} (2011), 185-216. \\
\ni [Lau] S. Lauritzen, {\sl Graphical models}.  Oxford University Press, 1996.\\
\ni [LeeJ] D. Lee and H. Joe, Multivariate extreme-value copulas with factor and tree dependence structures.  {\sl Extremes} {\bf 21} (2018), 147-176. \\
\ni [Loe] M. Lo\`eve, {\sl Probability theory}, Volumes I, II, 4th ed., Springer, 1977. \\
\ni [MaoH] T. Mao and T. Hu, Relations between the spectral measures and dependence of multivariate extreme value distributions.  {\sl Extremes} {\bf 18} (2015), 65-84. \\
\ni [MarP] M. B. Marcus and M. Pinsky, On the domain of attraction of $\exp(- e^{-x})$.  {\sl J. Math. Anal. Appl.} {\bf 28} (1968), 440-449.  \\
\ni [McNFE] A. J. McNeil, R. Frey and P. Embrechts, {\sl Quantitative risk management: Concepts, techniques, tools}.  Princeton University Press, 2005 (2nd ed. 2015). \\
\ni [McNN] A. J. McNeil and J. N\u eslehov\'a, Multivariate Archimedean copulas, $d$-monotone functions and ${\ell}_1$-norm symmetric distributions.  {\sl Ann. Stat.} {\bf 37} (2009), 3059-3097. \\
\ni [vMis] R. von Mises, La distribution de la plus grande de $n$ valeurs.   {\sl Revue Math\'ematique de l'Union Interbalkanique} {\bf 1} (1936), 141-160 (reprinted in {\sl Selected Papers II}, 271-294, AMS).   \\
\ni [Ost] A. J. Ostaszewski, Beurling regular variation, Bloom dichotomy, and
the Go\l\k{a}b-Schinzel functional equation.  {\sl Aequat. Math.} {\bf 89} (2015), 725-744. \\
\ni [PapT] I. Papastathopoulos and J. A. Tawn, Dependence properties of multivariate max-stable distributions.  {\sl J. Multivariate Anal.} {\bf 130} (2014), 134-140. \\
\ni [Pic1] J. Pickands III, Statistical inference using extreme order statistics.  {\sl Ann. Stat.} {\bf 3} (1975), 119-131. \\
\ni [Pic2] J. Pickands III, Multivariate extreme value distributions.  {\sl Bull. Inst. Int. Stat.} {\bf 49} (1981), 859-878. \\
\ni [Res1] S. I. Resnick, {\sl Extreme values, regular variation and point processes}.  Springer, 1987. \\
\ni [Res2] S. I. Resnick, {\sl Heavy-tail phenomena: Probabilistic and statistical modelling}.  Springer, 2007. \\
\ni [RooSW1] H. Rootz\'en, J. Segers and J. L. Wadsworth, Multivariate peaks
over thresholds models, {\sl Extremes} {\bf 21} (2018), 115-145. \\
\ni [RooSW2] H. Rootz\'en, J. Segers and J. L. Wadsworth, Multivariate generalized Pareto distributions: parametrizations, representations, and properties. {\sl J. Multivariate Analysis} {\bf 165} (2018), 117–-131. \\
\ni [RooT] H. Rootz\'en and N. Tajvidi, Multivariate generalized Pareto distributions.  {\sl Bernoulli} {\bf 12} (2006), 917-930. \\
\ni [Schl] M. Schlather, Models for stationary max-stable random fields.  {\sl Extremes} {\bf 5} (2002), 33-44. \\
\ni [SchmS] R. Schmidt and U. Stadtm\"uller, Non-parametric estimation of tail dependence.  {\sl Scand. J. Stat.} {\bf 33} (2006), 307-335. \\
\ni [Seg] J. Segers, Non-parametric inference for bivariate extreme-value copulas, {\sl Topics in Extreme Values} (ed. M. Ahsanullah and S.N.U.A. Kirmani), Nova Science, New York (2007), 181–-203. \\
\ni [ShaW] P. Sharkey and H. C. Winter, A Bayesian spatial hierarchical model for extreme precipitation in Great Britain. arXiv:1710.02091. \\
\ni [Sib] M. Sibuya, Bivariate extreme statistics.  {\sl Ann. Inst. Stat. Math.} {\bf 11} (1960), 195-210. \\
\ni [Skl] A. Sklar, Fonctions de r\'epartition \`a $n$ dimensions et leurs marges.  {\sl Publ. Inst. Stat. U. Paris} {\bf 8} (1959), 229-231. \\
\ni [Smi1] R. L. Smith, Max-stable processes and spatial extremes.  Unpublished manuscript, U. Surrey, 1990. \\
\ni [Smi2] R. L. Smith, Threshold methods for sample extremes. P.621-638 in {\sl Statistical extremes and aplications} (ed. J. Tiago de Oliveira), Reidel, 1984. \\

\bigskip

\bigskip

\bigskip

\bigskip

\ni N. H. Bingham, Mathematics Department, Imperial College, London SW7 2AZ; n.bingham@ic.ac.uk \\
\ni A. J. Ostaszewski, Mathematics Department, London School of Economics, Houghton Street, London WC2A 2AE; A.J.Ostaszewski@lse.ac.uk

\end{document}